\documentclass[a4paper,12pt]{article}
\usepackage{a4wide}
\usepackage{amsmath}
\usepackage{amssymb}
\usepackage{amsthm}
\usepackage{latexsym}
\usepackage{graphicx}
\usepackage[english]{babel}
\usepackage{makeidx}

\newtheorem{prop}[subsection]{Proposition}

\newtheorem{teor}[subsection]{Theorem}
\newtheorem{lema}[subsection]{Lemma}
\newtheorem{cor} [subsection]{Corollary}

\begin{document}

\selectlanguage{english}
\frenchspacing

\large
\begin{center}
\textbf{Some remarks on Borel type ideals.}

Mircea Cimpoea\c s
\end{center}

\normalsize

\footnotetext[1]{This paper was supported by the CEEX Program of the Romanian
Ministry of Education and Research, Contract CEX05-D11-11/2005 and by PN-II-...(to be completed after 15 September)}

\begin{abstract}
We give new equivalent characterizations for ideals of Borel type. Also, we prove that the 
regularity of a product of ideals of Borel type is bounded by the sum of the regularities of those ideals.

\vspace{5 pt} \noindent \textbf{Keywords:} Borel type ideals, Mumford-Castelnuovo regularity.

\vspace{5 pt} \noindent \textbf{2000 Mathematics Subject
Classification:}Primary: 13P10, Secondary: 13E10.
\end{abstract}

\section*{Introduction.}

Let $K$ be an infinite field, and let $S=K[x_1,\ldots,x_n],n\geq 2$ the polynomial ring over $K$.
Bayer and Stillman \cite{BS} note that Borel fixed ideals $I\subset S$ satisfy the following property:
\[(*)\;\;\;\;(I:x_j^\infty)=(I:(x_1,\ldots,x_j)^\infty)\;for\; all\; j=1,\ldots,n.\] Herzog, Popescu and Vladoiu \cite{hpv} define a monomial ideal $I$ to be of \emph{Borel type} if it satisfies $(*)$. We mention that this concept appears also in \cite[Definition 1.3]{CS} as the so called {\em weakly stable ideal}. Herzog, Popescu and Vladoiu proved in \cite{hpv} that $I$ is of Borel type, if and only if for any monomial $u\in I$ and for any $1\leq j<i \leq n$  with $x_i|u$, there exists an integer $t>0$ such that $x_j^{t}u/x_i^{\nu_i(u)}\in I$, where $\nu_i(u)$ is the exponent of $x_i$ in $u$. This allows us to prove that the property of an ideal to be of Borel type is preserved for several operations, like sum, intersection, product, colon, see Proposition $1.1$.

Let $G(I)$ be the minimal set of monomial generators of $I$ and $deg(I)$ the highest degree of a monomial of $G(I)$. Given a monomial $u \in S$ we set $m(u) = max\{i :\; x_i|u\}$ and \linebreak $m(I) = max_{u\in G(I)} m(u)$. If $\beta_{ij}(I)$ are the graded Betti numbers of $I$ then, the regularity of $I$ is given by 
$reg(I) = \max\{j-i :\; \beta_{ij}(I)\neq 0\}$.
S.Ahmad and I.Anwar presented in \cite{saf} a characterization for ideals of Borel type and proved that if $I\subset S$ is an ideal of Borel type, then $reg(I) \leq m(I)(deg(I)-1)+1$, see
\cite[Cor 2.4]{saf}. On the other hand, we proved in \cite{mir} that $reg(I)=min\{e\geq deg(I):\; I_{\geq e}$ is stable $\}$, see \cite[Cor 8]{mir}. We extend this result, showing that $I$ is an ideal of Borel type if and only if there exists some integer $e\geq deg(I)$ such that $I_{\geq e}$ is stable, see Theorem $1.5$.

Using the above results, we prove that the regularity of a product of ideals of Borel type is bounded by the sum of the regularity of each ideal which appear in that product, see Theorem $1.7$. As a consequence, if $I$ is an ideal of Borel type then $reg(I^k)\leq k\cdot reg(I)$ for any positive integer $k$, see Corollary $1.8$. In the second section of our paper, we give an explicit description for the ideals of Borel type, see Theorem $2.1$.

\newpage
\section{Remarks on Borel type ideals.}

\begin{prop}
(1) If $I,J\subset S$ are two ideals of Borel type then $I+J$, $I\cap J$ and $I\cdot J$ are also ideals of Borel type.

(2) If $I\subset S$ is an ideal of Borel type and $J\subset S$ is an arbitrary monomial ideal, then $(I:J)$ is an ideal
    of Borel type.
\end{prop}

\begin{proof}
(1) Since a monomial of $I+J$ is either in $I$, either in $J$ it follows immediately that $I+J$ is of Borel type.
    A similar argument
    holds for $I\cap J$. Now, let $u\in I\cdot J$ be a monomial. It follows that $u=v\cdot w$, where $v\in I$ and $w\in J$
    are monomials. Let $1\leq i\leq n$ such that $x_i|u$ and let $1\leq j<i$. If $x_i|v$, since $I$ is of Borel type, 
    then there exists some integer $t_1>0$ such that $x_j^{t_1}\cdot v/x_i^{\nu_i(v)}\in I$. Otherwise, if $x_i$ does not
  divide $v$ we put $t_1=0$. Analogously, there exists some integer $t_2\geq 0$ such that $x_j^{t_2}\cdot
   w/x_i^{\nu_i(w)}\in J$.
    It follows that $x_j^{t_1+t_2}u/x_i^{\nu_i(u)}\in I\cdot J$, therefore $I\cdot J$ is of Borel type.

(2) Suppose $J=(v_1,\ldots,v_m)$, where $v_i$ are monomials. Since $(I:J)=\cap_{k=1}^{m} (I:v_k)$ and the intersection
    of Borel type ideals is still of Borel type, we can assume $m=1$. Denote $v_1:=v$. Let $u\in (I:v)$. We have
    $u\cdot v\in I$. Let $1\leq i\leq n$ such that $x_i|u$ and let $1\leq j<i$. Since $I$ is of Borel type, there
    exists some $t\geq 0$ such that $x_j^t u\cdot v/x_i^{\nu_i(uv)}\in I$. In particular, multiplying by
    $x_i^{\nu_i(v)}$, it follows that $v\cdot(x_j^t u/x_i^{\nu_i(u)})\in I$ and thus $x_j^t u/x_i^{\nu_i(u)}\in (I:v)$.
    In conclusion, $(I:v)$ is of Borel type as required.
\end{proof}

As a consequence, we give a new proof of the following well known fact, see
for instance \cite[Proposition 5.2]{hp}.

\begin{cor}
If $I\subset S$ is an ideal of Borel type, then any associated prime $P\in Ass(I)$ has the form $I=(x_1,\ldots,x_r)$
for some $1\leq r\leq n$.
\end{cor}

\begin{proof}
Any prime ideal $P\in Ass(I)$ can be written as $P=(I:u)$ for some monomial $u\in S$. Since $I$ is of Borel type, 
by Proposition $1.1$ we get $P$ of Borel type. On the other hand, $P$ is a prime monomial ideal, so it is generated by
variables. Combining these facts, we are done.
\end{proof}

\begin{lema}
Let $I\subset S$ be a monomial ideal. The following are equivalent:

(1) $I$ is of Borel type,

(2) For any monomial $u\in I$ and
for any integer $1\leq j<m(u)$, there exists an integer $t>0$ such that $x_j^t u/x_{m(u)}^{\nu_{m(u)}(u)} \in I$.
\end{lema}

\begin{proof}
The first implication is trivial. For the converse implication, let $u\in I$ be a monomial and let $j<i$ such that $x_i|u$. Suppose $u=x_1^{a_1}\cdots x_i^{a_i}\cdots x_{q}^{a_q}$, where $q=m(u)$ and $a_1,\ldots,a_q$ are some nonnegative integers. By our assumption, there exists an integer $t_q>0$ such that $x_j^{t_q}u/x_{q}^{a_q} = x_j^{t_q}x_1^{a_1}\cdots x_{q-1}^{a_{q-1}} \in I$. Recursively, we can find some positive integers $t_{q-1},\ldots,t_i$ such that $x_j^{t_q+t_{q-1}+\cdots+t_i} x_1^{a_1}\cdots x_{i-1}^{a_{i-1}}\in I$. We denote $t=t_q+t_{q-1}+\cdots+t_i$. It follows that,
\[ x_j^tu/x_i^{a_i}:= (x_j^{t_q+t_{q-1}+\cdots+t_i} x_1^{a_1}\cdots x_{i-1}^{a_{i-1}})\cdot(x_{i+1}^{a_{i+1}}\cdots
x_q^{a_q})\in I, \] thus $I$ is of Borel type.
\end{proof}

We recall the following result of Eisenbud-Reeves-Totaro, see \cite[Proposition 12]{ert}.

\begin{prop}
Let $I\subset S$ be a monomial ideal and let $e\geq deg(I)$ be an integer. If 
$I_{\geq e}$ is stable then $reg(I)\leq e$.
\end{prop}

In the following, we will give a theorem of characterization for ideals of Borel type.

\begin{teor}
Let $I\subset S$ be a monomial ideal. The following are equivalent:

(1) $I$ is an ideal of Borel type.

(2) $I_{\geq reg(I)}$ is stable.

(3) There exists an integer $e\geq deg(I)$ such that $I_{\geq e}$ is stable.

\noindent
If one of the above conditions holds, $reg(I)=min\{e|\;e\geq deg(I)\;, I_{\geq e}$ is stable $\}$.
\end{teor}

\begin{proof}
$(1)\Rightarrow (2)$ follows from \cite[Theorem 6]{mir}. $(2)\Rightarrow (3)$ is obvious. In order to prove $(3)\Rightarrow (1)$, we choose a monomial $u\in I$. We denote $q=m(u)$ and thus $u=x_1^{a_1}\cdots x_q^{a_q}$ for some nonnegative integers $a_i$ with $a_q>0$. Let $j<q$ be an integer. According with Lemma $1.3$, we must show that there exists some integer $t>0$ such that $x_j^t\cdot u/x_q^{a_q}\in I$. 

We choose an integer $t'>0$ such that $deg(x_j^{t'}u)\geq e$
and therefore $x_j^{t'}\cdot u \in I_{\geq e}$. Since $I_{\geq e}$ is stable, it follows that $x_j^{t'+1}\cdot u/x_q\in I$, $x_j^{t'+2}\cdot u/x_q^2\in I, \ldots$ and finally, $x_j^{t'+a_q}u/x_q^{a_q}\in I$. Denoting $t=t'+a_q$ we get the required conclusion.

The last assertion of the theorem follows immediately from Proposition $1.4$ and $(2)$.
See also \cite[Cor 8]{mir}.
\end{proof}

\begin{lema}
Let $I,J\subset S$ be two monomial ideals and let $e\geq deg(I)$ and $f\geq deg(J)$ be two integers such that 
$I_{\geq e}$ and $J_{\geq f}$ are stable. Then $(I\cdot J)_{\geq e+f}$ is stable.
\end{lema}

\begin{proof}
Let $u\in (I\cdot J)_{\geq e+f}$ be a monomial. It follows that $u=v\cdot w$ for some monomials $v\in I$ and $w\in J$.
Since $deg(u)\geq e+f$, we can easily choose the monomials $v$ and $w$ such that $v\in I_{\geq e}$ and $w\in J_{\geq e}$. 

Now, let $j<m(u)$ be an integer and suppose $m(u)=m(v)$. Then $x_j u/x_{m(u)} = (x_j v/x_{m(v)})\cdot w \in I\cdot J$,
because $x_j v/x_{m(v)}\in I$ since $I_{\geq e}$ is stable. Analogously, if $m(u)=m(w)$, then
$x_j u/x_{m(u)} = v\cdot (x_j w/x_{m(w)}) \in I\cdot J$. Therefore, $(I\cdot J)_{\geq e+f}$ is stable.
\end{proof}

\begin{teor}
Let $I,J\subset S$ be two monomial ideals of Borel type. Then
\[ reg(I\cdot J) \leq reg(I) + reg(J).\]
\end{teor}

\begin{proof}
Since $I$ and $J$ are ideals of Borel type, if we denote $e:=reg(I)$ and $f=reg(J)$, by Theorem $1.5(2)$, or 
\cite[Theorem 6]{mir}, we get that $I_{\geq e}$ and $J_{\geq f}$ are stable. Using the previous lemma, it follows that $(I\cdot J)_{\geq e+f}$ is stable, therefore, using Proposition $1.3$ we get $reg(I\cdot J)\leq e+f$ as required.
\end{proof}

\begin{cor}
If $I\subset S$ is an ideal of Borel type, then $reg(I^k)\leq k\cdot reg(I)$.
\end{cor}

Note that, there are other large classes of graded ideals which have this property, see for instance \cite{CH}, but on the other hand, Sturmfels provided an example of a graded ideal $I\subset S$ with $reg(I^2)>2\cdot reg(I)$ in \cite{S}.

\section{An explicit description of Borel type ideals.}

\begin{teor}
Let $I\subset S$ be an ideal of Borel type. Then, there exists an integer $1\leq q\leq n$, 
some nonnegative integers $r_0,r_1,\ldots,r_{n-q}$, and some monomials
$v_{0j}\in K[x_{1},\ldots,x_{q}]$ for $j\in \{1,\ldots,r_0\}$ and $v_{ij}\in K[x_{1},\ldots,x_{q+i-1}]$ for any $1\leq i\leq n-q$ and $j\in \{1,\ldots,r_i\}$, such that the set of minimal monomial generators of $I$, $G(I)$ is:
\[ \{ x_1^{a_1},\ldots,x_q^{a_q},v_{01},\ldots,v_{0r_0},v_{11}x_{q+1}^{a_{11}},\ldots,
  v_{1r_1}x_{q+1}^{a_{1r_1}},\ldots, v_{n-q,1}x_{n}^{a_{n-q,1}},\ldots, v_{n-q,r_{n-q}}x_n^{a_{n-q,r_{n-q}}} \}, \]
where $a_j$ for $j=\overline{1,q}$ and $a_{ij}$ for $i=\overline{1,n-q}$, $j=\overline{1,r_{i}}$ are some positive integers. Moreover,if $r_i>0$ for some $i\geq 1$ it follows that $r_{i+1}>0$. Also, for any $2\leq i\leq n-q$ and $1\leq j\leq r_i$ there exists some $1\leq k\leq r_{i-1}$ such that $v_{i-1,k}|v_{ij}$.

Conversely, any ideal monomial ideal $I\subset S$ with $G(I)$ satisfying the conditions above, is an ideal of Borel type.
\end{teor}

\begin{proof}
Let $q:=max\{j|\; x_j\in \sqrt{I}\}$. It follows that there exists a positive integer $k$ such that
$x_q^{k}\in I$. We choose $a_q:=min\{k:\;x_q^k\in I\}$. Therefore, we have $x_q^{a_q}\in G(I)$.
Since $I$ is of Borel type, it follows that there exists some positive integers $a_1,\ldots,a_{q-1}$
such that $x_j^{a_j}\in I$ for any $1\leq j\leq q-1$. As above, we can assume $x_j^{a_j}\in G(I)$ for any $1\leq j\leq q-1$. We denote by $v_{01},\ldots,v_{0r_0}$ those monomials from $G(I)$ which are in $K[x_1,\ldots,x_q]$, but are not in the set $\{x_1^{a_1},\ldots,x_q^{a_q}\}$, if any.

Let $1\leq i\leq n-q$ be an integer. Suppose $\{u\in G(I)|\; m(u)=q+i\} = \{w_{i1},\ldots,w_{ir_i}\}$. It follows
that $w_{ij}=v_{ij}x_{i+q}^{a_{ij}}$, where $v_{ij}\in K[x_1,\ldots,x_{i+q-1}]$ and $a_{ij}>0$.
Assume $r_i=0$ for some $0\leq i\leq n-q-1$. Suppose $r_{i+1}>0$ and therefore there exists a monomial $v_{i+1,1}\in K[x_1,\ldots,x_{q+i}]$ and a positive integer $a_{i+1,1}$ such that $v_{i+1,1}\cdot x_{i+q+1}^{a_{i+1,1}}\in G(I)$. Since
$I$ is of Borel type, it follows that there exists a positive integer $t>0$ such that $v_{i+1,1}x_{q+i}^{t}\in I$. We can write $v_{i+1,1}x_{q+i}^{t} = v\cdot x_{q+i}^{t'}$ such that $m(v)<k$, where $v$ is a monomial. Now, since $v\cdot x_{q+i}^{t'}\in I$ and $r_{q+i}=0$, it follows that $v$ is a multiple of a minimal generator of $I$ and moreover $v_{i+1,1}$ is multiple of a minimal generator of $I$. This is a contradiction with the fact that $v_{i+1,1}x_{q+i+1}^{a_{i+1,1}}$ is a minimal generator of $I$.

Suppose $r_i>0$ for some $i\geq 2$ and let $1\leq j\leq r_i$ be an integer. Since $v_{ij}x_{q+i}^{a_{ij}}\in I$ and
$I$ is of Borel type, it follows that $v_{ij}x_{q+i-1}^{t}\in I$ for some $t>0$. It follows that there exists $u\in G(I)$ with $u|v_{ij}x_{q+i-1}^{t}$. Obviously, we cannot have $u\in K[x_1,\ldots,x_{q+i-2}]$ because, in such case,
it follows that $u|v_{ij}$ and we contradict the minimality of $v_{ij}x_{q+i}^{a_{ij}}$ as a generator of $I$. Therefore, $u=v_{i-1,k}x_{q+i-1}^{a_{i-1,k}}$ for some $k\leq r_{i-1}$ and thus thus $v_{i-1,k}|v_{ij}$.

For the converse, we consider an ideal $I$ with $G(I)$ fulfilling the conditions of the theorem 
and we must show that $I$ is of Borel type. Let $u\in I$ be a monomial. Let $s<m(u)$ be an integer. If $s\leq q$ we choose $t:=a_s$ and we get that $x_s^{a_s}u/x_{m(u)}^{\nu_{m(u)}(u)}\in I$. Now, we assume $s\geq q+1$ and we denote $j:=s-q$. Also, we can assume that $u$ is a minimal generator of $I$ and thus $u = v_{im}x_{q+i}^{a_{im}}$ for some $1\leq i\leq n-q$ and $1\leq m\leq r_i$. Note that $m(u)=q+i$ and moreover, $\nu_{q+i}(u)=a_{im}$. By Lemma $1.3$, 
in order to show that $I$ is of Borel type, it is enough to prove that there exists some integer $t>0$ such that $x_{q+j}^t \cdot v_{im} \in I$. Indeed, by our assumption, it follows that there exists some 
$k\leq r_j$ such that $v_{jk}|v_{im}$. But, $v_{jk}x_{q+j}^{a_{jk}}\in G(I)$, therefore, if we choose $t\geq a_{jk}$,
it follows that $x_{q+j}^t \cdot v_{im} \in I$, as required. 
\end{proof}

\vspace{2mm} \noindent {\footnotesize
\begin{minipage}[b]{15cm}
 Mircea Cimpoea\c s, Institute of Mathematics of the Romanian Academy, Bucharest, Romania\\
 E-mail: mircea.cimpoeas@imar.ro
\end{minipage}}
\end{document}